\documentclass[11pt,a4paper]{amsart}

\usepackage{amssymb,amsmath,graphics,verbatim}
\usepackage{latexsym}
\usepackage{eucal}
\usepackage{mathrsfs}
\usepackage{a4wide}

\usepackage{color}

\newtheorem{theorem}{Theorem}

\newtheorem{lemma}[theorem]{Lemma}
\newtheorem{prop}[theorem]{Proposition}

\theoremstyle{remark}

\newcommand{\E}{\mathbb{E}}
\newcommand{\Prob}{\mathbb{P}}
\newcommand{\R}{\mathbb{R}}
\newcommand{\N}{\mathbb{N}}
\newcommand{\C}{\mathbb{C}}
\newcommand{\T}{\mathbb{T}^2}

\newcommand{\Z}{\mathbb{Z}}

\newcommand{\lt}{\mathscr{L}}

\newcommand{\lp}{\mathcal{L}}
\newcommand{\Aniso}{\mathcal{H}_{r,M}}
\newcommand{\ft}{\mathcal{F}_\theta}
\newcommand{\ftp}{\mathcal{F}^+_\theta}

\begin{document}
\bibliographystyle{plain}

\title{On the rate of mixing of Circle extensions of Anosov maps.}
\date{\today}
\author{Fr\'ed\'eric Naud}
\address{Laboratoire de Math\'ematiques d'Avignon \\
Campus Jean-Henri Fabre, 301 rue Baruch de Spinoza\\
84916 Avignon Cedex 9. }
\email{frederic.naud@univ-avignon.fr}

\begin{abstract}
Let $A:\T\rightarrow \T$ be an Anosov diffeomorphism. Circle extensions $\widehat{A}$ are a rich family of non-uniformly hyperbolic diffeomorphisms
living on $\T\times S^1$ for which the rate of mixing is conjectured to be generically exponential. 
In this paper, we investigate the possible rates of exponential mixing by exhibiting
some {\it explicit lower bounds} on the decay rate by {\it Fourier analytic} and {\it probabilistic} techniques. The rates obtained are related to the topological pressure of two times the unstable jacobian.
\end{abstract}

\maketitle

\section{Introduction}
 Let $\T=\R^2\slash \Z^2$ be the usual flat torus. And consider $A:\T\rightarrow \T$ an Anosov diffeomorphism, which will assumed
 to be topologically mixing in the sequel.
 Let $\tau:\T\rightarrow \R$ be a smooth map and let $S^1:=\R\slash \Z$ be the circle. Then one can define
 {\it an $S^1$-extension of $A$}, denoted by $\widehat{A}_\tau:\T\times S^1\rightarrow \T\times S^1$  by setting
 $$\widehat{A}_\tau(x,\omega):=(Ax,\tau(x)+\omega),$$
 all coordinates being understood mod $2\pi$. 
 These maps $\widehat{A}_\tau$ are the simplest prototype of {\it partially} hyperbolic systems, for which the neutral
 direction forms a trivial bundle in the tangent bundle. The qualitative ergodic theory of these maps is well established, and most questions of ergodic stability are settled in the work of Brin \cite{Brin1} and Burns-Wilkinson \cite{BurnsWilkinson}. 
 However, when it comes to quantitative ergodic theory,
 very few results are known. Let $\mu_{srb}$ be the {\it Sinai-Ruelle-Bowen} $A$-invariant probability measure, which
 can be characterized as the unique {\it physical measure}, for which Birkhoff averages converge Lebesgue-almost surely
 to the spatial average. A natural $\widehat{A}$-invariant extension of $\mu_{srb}$ to $\T\times S^1$ can be defined by
 $$\int F(x,\omega) d\widehat{\mu}_{srb}(x,\omega):=\int_{\T}\int_{S^1} F(x,\omega)d\mu_{srb}(x)d\omega.$$
 From the pioneering work of Dolgopyat \cite{Dolgopyat}, it follows 
 that for generic $\tau$, the map $\widehat{A}$ has {\it rapid decay of correlations} for all $C^\infty$ observables, i.e. for all $\varphi,\psi$ $C^\infty$ on $\T\times S^1$, we have as $N\rightarrow +\infty$ and all $k\geq 1$, 
 $$C_{\varphi,\psi}(N):=\int \left (\varphi \circ \widehat{A}^N \right) \psi d\widehat{\mu}_{srb} =
 \int \varphi d\widehat{\mu}_{srb}\int \psi d\widehat{\mu}_{srb}+O_{\varphi,\psi,k}(N^{-k}).$$
 It is natural to expect that {\it exponential mixing} is also typical, but it is still an open question in the context of extensions. On the other hand, for Anosov Flows,
 a recent preprint of Tsujii \cite{Tsujii} shows that generic volume preserving $3$-dimensional Anosov flows are exponentially mixing. See also \cite{BW1} for new results in higher dimensions. We don't know if Tsujii's recent technique can be used to prove exponential mixing in our context, and this should be pursued elsewhere.

 One natural question raised by our current knowledge is what's the typical rate of mixing when observables are very regular ?
 Could we get super exponential mixing as in the case of linear Anosov diffeomorphisms of $\T$ ? How does
 the instability of the system (Lyapunov exponents) affect this rate of mixing ? To formulate our main result,
 we recall that the {\it unstable jacobian} is defined by
 $$J^u(x)=\vert \det ( D_x A\vert_{E^u_x})\vert=\vert D_x A \vert_{E^u_x},$$
 where $E^u_x\subset T_x \T$ is the unstable direction at $x$ of $A$, which is $D_x A$-invariant. This is at least a H\"older continuous function on $\T$. Given a H\"older function $\varphi$ on $\T$ one can define the topological pressure $P(\varphi)$ by taking the supremum
 $$P(\varphi):=\sup_{\mu\  A-\mathrm{inv}}\left ( h_\mu(A)+\int_{\T} \varphi d\mu \right),$$
 where the sup is taken over all invariant probability measures, and $h_\mu(A)$ is the measure theoretic entropy of $\mu$.
 Our main result is the following.
 \begin{theorem}
 \label{main}
  Assume that $A$ is a $C^1$-small enough, {\bf volume preserving} real analytic perturbation of a linear Anosov map. 
 \begin{itemize}
 \item Then for all {\bf $\tau:\T\rightarrow \R$ real analytic}, for all $\epsilon>0$, one can find real analytic observables $\varphi,\psi$  with $\int \varphi=\int \psi=0$ such that 
  $$\limsup_{N\rightarrow +\infty}  \vert C_{\varphi,\psi}(N) \vert ^{\frac{1}{N}}\geq e^{\frac{5}{2} P(-2\log J^u)-\epsilon}.$$
  \item For all $N$ large enough, almost surely for all $\tau(z)=P_N(z)\in \mathcal{P}_N$ random trigonometric polynomial, the extended map $\widehat{A}_\tau$ is rapidly mixing for the SRB-measure on $\T\times S^1$.
  \item For all $\epsilon>0$, with positive probability for $\tau \in \mathcal{P}_N$,
  one can find real analytic observables $\varphi,\psi$  with $\int \varphi=\int \psi=0$ and
  $$\limsup_{N\rightarrow +\infty}  \vert C_{\varphi,\psi}(N) \vert ^{\frac{1}{N}}\geq e^{\frac{1}{2} P(-2\log J^u)-\epsilon}.$$
  \end{itemize}
 \end{theorem}
 The set of random trigonometric polynomials of degree $N$, denoted by $\mathcal{P}_N$, is defined in $\S 4$ and is just the obvious
 guess: independent Gaussian combinations of Laplace eigenfunctions on $\T$.
 The above theorem shows in particular that the rate of mixing, unlike in the uniformly hyperbolic case, can never be super exponential. 
 This fact was already pointed out for suspensions of analytic expanding maps by the author in \cite{NaudHP}, with a less precise lower bound
 involving entropy rather than pressure. Note that the first statement is unconditional and holds for all choice of $\tau$ but there is a loss in the lower bound. However, we are also able to show that for "many" choices of $\tau$ among the set of random polynomials $\mathcal{P}_N$, the rate of mixing is bigger than $$e^{\frac{1}{2} P(-2\log J^u)-\epsilon},$$ 
 which we believe is the optimal lower bound. One of the motivations for this type of quantitative lower bounds is that it shows that
 when unstable Lyapunov exponents are "small", i.e. close to $1$, then the rate of mixing is arbitrarily (exponentially) slow. 
 
 \bigskip
 \noindent {\bf An example: Arnold's cat map.} This is the standard Anosov diffeomorphism on $\T$ induced by the action of the $SL_2(\Z)$
 matrix
 $$M=\left ( \begin{array}{cc} 2&1\\1&1 \end{array}    \right).$$
 The eigenvalues are $\lambda^{\pm}=\frac{3\pm \sqrt{5}}{2}$ which implies that the topological entropy $h(M)$ of the cat map is exactly
 $$h(M)=\log\left(\frac{3+\sqrt{5}}{2}\right).$$
The topological pressure $P(-2\log J^u)$ is easily computed as 
$$P(-2\log J^u)=h(M)-2\log(\lambda^+)=-\log\left(\frac{3+\sqrt{5}}{2}\right), $$ 
which yields 
$$\exp\left(\frac{1}{2} P(-2\log J^u)\right)=\sqrt{\frac{2}{3+\sqrt{5}}}\approx 0,618033988.$$
Theorem \ref{main} tells us that while the map $M$ itself mixes at super exponential decay rate for all analytic observables (for an elementary
proof of that fact, see \cite{BaladiBook1}, chapter 4), there exists (at least rapidly mixing) extensions $\widehat{M}_\tau$ and analytic  observables whose rate of mixing is not faster than $(0,618)^N$.

\bigskip
The paper is organized as follows. In the next section, we show how the lower bounds on correlation functions can be derived from
a statement on the spectrum of certain twisted transfer operators $\lt_q$ that depend on a frequency parameter $q\in \Z$. These operators
act naturally on an anisotropic function space defined by Faure-Roy in \cite{FaureRoy}. All the material and a priori estimates regarding
these spaces is gathered in the last section $\S 5$. In $\S 3$, we prove the first part of the main spectral estimate via a technique based
on "frequency averaging", i.e. we prove certain bounds by summing smoothly over $q$ and eventually recover some pointwise bounds.
In $\S 4$, we use a different averaging technique with a more probabilistic flavour: we consider some random "roof functions" $\tau$ (given by a random ensemble of trigonometric polynomials of degree $N$) and show that one can prove a lower bound on the {\it expectation of the spectral radius} 
$\rho(\lt_q)$. This in turn shows the existence of a set of functions $\tau$ with an improved lower bound on the spectral radius. We also prove, using mostly
the old technology of subshifts of finite type, that provided the degree $N$ is taken large enough, exponential mixing occurs with probability $1$ in 
$\mathcal{P}_N$.
The last section $\S5$ is devoted to properties of the Anisotropic function space $\Aniso$ and we rely on the existing work \cite{AlexAdam,FaureRoy} and provide proofs of some spectral upper bounds that are necessary for our purpose.

\section{Function space and reduction to a spectral problem}
The main result (Theorem \ref{main}) follows from a statement on the spectrum of certain "twisted" transfer operators. First, observe
that given an observable $F(x,\theta)$ defined on $\T\times S^1$ of the form ($q\in \Z$)
$$F(x,\theta)=f(x)e^{iq2\pi\theta},$$
then we have
$$F\circ \widehat{A}(x,\theta)=  \left(e^{i2\pi q\tau(x)}f\circ A(x)\right) e^{iq2\pi\theta}.$$ 
this leads naturally to study the following "twisted" koopman operators $\lt_q$ acting by
$$\lt_q(f)(x):=e^{iq 2\pi\tau(x)}f\circ A(x).$$

The analysis of $\lt_q$ will depend crucially on a good choice of function space. We will be working in the real-analytic category
and we describe below the functional analytic set-up. Let $M\in \mathrm{SL}_2(\Z)$ be a hyperbolic matrix so that its induced action
on $\T$ is an Anosov map. Let $r>0$ be a parameter. We recall that a trigonometric polynomial $P(x)$ on $\T$ is simply an expression
of the type
$$P(x)=\sum_{\vert \alpha \vert \leq N} a_\alpha e^{i2\pi \alpha .x},$$
where $\alpha=(\alpha_1,\alpha_2)\in \Z^2$, $\vert \alpha \vert=\vert \alpha_1\vert +\vert \alpha_2\vert$ and
$\alpha .x=\alpha_1x_1+\alpha_2 x_2 $. Trigonometric polynomials are obviously real-analytic on $\T$ and extend holomorphically to $\C^2$.
\begin{theorem}
\label{blackbox} 
There exists a family of Hilbert spaces $\mathcal{H}_{r,M}$ which contain densely all trigonometric polynomials on $\T$, such that we have:
\begin{enumerate} 
\item For all $\tau$ real analytic on $\T$, for all $C^1$-small enough real analytic perturbation $A$ of $M$, one can find $r>0$ such that
$\lt_q:\mathcal{H}_{r,M}\rightarrow \mathcal{H}_{r,M}$ acts as a bounded compact trace class operator.
\item For all $q$, the spectral radius $\rho(q)$ of $\lt_q$ is smaller than $1$.
\item Moreover there exist constants $C,\beta>0$, and $r'>r$, independent of $q$ such that the eigenvalue sequence $\lambda_k(\lt_q)$ satisfies the bound
$$\vert \lambda_k(\lt_q)\vert \leq Ce^{\vert q \vert \Vert \tau \Vert_{r',\infty}} e^{-\beta \sqrt{k}}.$$
\item The lebesgue measure $dm$ on the torus extends as a continuous linear functional $L_{m}:\mathcal{H}_{r,M}\rightarrow \mathcal{H}_{r,M} $ and has "full support" in the following generalized sense. Given  $\varphi \in \mathcal{H}_{r,M}$ with $\varphi\neq 0$, one can find a trigonometric polynomial $g$ such that \footnote{The fact that given $g\in \Aniso$, for all trigonometric polynomial $\varphi$ the product $g\varphi$ belongs to
$\Aniso$ will be clarified in $\S5.1$.}
$$L_{m}(g\varphi)\neq 0 .$$
\item For all $n$, we have the trace formula
$$\mathrm{Tr}( \lt_q^n)=\sum_{A^nx=x} \frac{e^{2i\pi q \tau^{(n)}(x)}}{\vert \det(I-D_xA^n)\vert},$$
where the sum runs over all periodic points of period $n$ of the map $A^n:\T\rightarrow \T$, and 
$\tau^{(n)}(x)=\tau(x)+\tau(Ax)+\ldots\tau(A^{n-1}x)$.
\end{enumerate}
\end{theorem}
Here the norm $\Vert \tau \Vert_{r',\infty}$ refers to the sup norm of $\tau$ in a {\it complexified neighbourhood} of the torus $\T$. More precisely if $f:\R^2\rightarrow \C$ extends holomorphically in small complex neighbourhood of the type
$\R^2+i[-r',+r']^2\subset \C^2$, we set
$$\Vert f \Vert_{r',\infty}=\sup_{x\in \T,y\in [-r',+r']^2} \vert f(x+iy) \vert.$$
The existence of such function spaces "adapted to the hyperbolic dynamics" follow in the analytic category from the work of Faure-Roy \cite{FaureRoy}. More recently, these function spaces have been used to study the Ruelle spectrum of Anosov maps by Adam 
\cite{AlexAdam}, and also Bandtlow-Just-Slipantschuk \cite{BJS2}. We will provide more details for the construction of these spaces
later on but roughly speaking, they are designed in Fourier coordinates to impose analyticity in the stable direction (exponential decay
of Fourier coefficients) while irregularity is allowed in the unstable direction (exponential growth at most of Fourier modes).
Part $1)$ and $5)$ follow readily from the above mentioned papers.
On the other hand parts $2)$ and $3),4)$ of the above theorem will require an extra amount of work which is the purpose of the last section of the paper. However, we can use Theorem \ref{blackbox}  as a "blackbox" to prove our main result which will follow from the spectral
statement below. 
\begin{prop}
\label{mainprop}
Under the above assumptions and notations, the following holds.
\begin{enumerate}
\item For all $\epsilon >0$, there exist infinitely many $q\in \Z$ such that 
$\lt_q:\mathcal{H}_{r,M}\rightarrow \mathcal{H}_{r,M}$ has an eigenvalue $\lambda(q)$ with
$$\vert  \lambda(q)\vert \geq e^{\frac{5}{2} P(-2\log J^u)-\epsilon}.$$
\item Furthermore, for all $\epsilon>0$, there exist trigonometric polynomials $\tau$, non cohomologous to constants, such that one can find $q\neq 0$ such that 
$\lt_{q,\tau}:\mathcal{H}_{r,M}\rightarrow \mathcal{H}_{r,M}$ has an eigenvalue $\lambda(q)$ with
$$\vert  \lambda(q)\vert \geq e^{\frac{1}{2} (P(-2\log J^u)-\epsilon)}.$$
\end{enumerate}
\end{prop}
Let us show how to derive Theorem \ref{main} from Proposition \ref{mainprop}. We fix $q\neq 0$ large enough such that either $1)$ or $2)$ from the above statement holds. We use observables of the form
$$\varphi(x,\theta)=f(x)e^{2i\pi q\theta};\ \psi(x,\theta)=g(x) e^{-2i\pi q\theta},$$
where $f,g$ will be analytic functions on the torus specified later on. Notice that we have 
$$\int \varphi d\widehat{\mu}_{srb}= \int \psi d\widehat{\mu}_{srb}=0.$$
Because $\lt_q:\mathcal{H}_{r,M}\rightarrow \mathcal{H}_{r,M}$ is compact, we can use holomorphic functional calculus to write for all $\epsilon,\rho>0$
$$\lt_q^N=\sum_{\vert \lambda_j(q) \vert >\rho} \lt_q^N \mathcal{P}_j +O((\rho+\epsilon)^N),$$
the error term being understood in the operator norm topology. Each $\mathcal{P}_j$ is a finite rank projector such that 
$$d_j:=\mathrm{ dim}(Im( \mathcal{P}_j))$$
equals the algebraic multiplicity\footnote{We cannot discard the possible presence of Jordan blocks.} 
of the eigenvalue $\lambda_j(q)$. For all $j\neq j'$ we have $\mathcal{P}_j\mathcal{P}_{j'}=0$, and 
$$\lt_q \vert_{ Im( \mathcal{P}_j)}=\lambda_j Id+\mathcal{N}_j, $$
where $\mathcal{N}_j^{d_j}=0$. Going back to the correlation function (from now on we assume that we are in the volume preserving case i.e. $d\mu_{srb}=dm$), we have
$$C_{\varphi,\psi}(N)=\int_{\T} \lt_q^N(f) g dm= \sum_{\vert \lambda_j(q) \vert >r} L_{m}(g\lt_q^N \mathcal{P}_j f) +O((\rho+\epsilon)^N)$$
$$=\sum_{\vert \lambda_j(q) \vert >\rho} \lambda_j^N Q_{j,f,g}(N) +O((\rho+\epsilon)^N),$$
where $Q_{j,f,g}(N)$ is a polynomial in $N$ with degree at most $d_j$, given by
$$Q_{j,f,g}(N)=L_{m}(\mathcal{P}_j f)+N\lambda_j^{-1}L_{m}(\mathcal{N}_j \mathcal{P}_j f)+\ldots+
\binom{N}{d_{j}-1} \lambda_j^{-d_j +1}L_{m}(g \mathcal{N}_j^{d_j-1} \mathcal{P}_j f).$$
Now set $\widetilde{\rho}:=\max_j \vert \lambda_j(q) \vert$. We set for all $j$ such that $\vert \lambda_j(q)\vert=\widetilde{\rho}$, $\lambda_j=\widetilde{\rho} e^{i \theta_j}$. We can then write 
$$C_{\varphi,\psi}(N)=\widetilde{\rho}^N\sum_{\vert \lambda_j(q) \vert =\widetilde{\rho}} e^{iN\theta_j} Q_{j,f,g}(N) +O((\rho)^N) $$
for some $\rho< \widetilde{\rho}$. Let us set \footnote{Remark that if there is only one eigenvalue with $\vert \lambda_j \vert =\widetilde{\rho}$, then the proof is much simpler, but we cannot exclude this case.}
$$S(N):= \sum_{\vert \lambda_j(q) \vert =\widetilde{\rho}} e^{iN\theta_j} Q_{j,f,g}(N).$$ 
To obtain a {\it lower bound } on the oscillating sum $\vert S(N)\vert$, we will use Dirichlet box principle, which for us is the following handy fact.

\begin{lemma}
 Let $\alpha_1,\ldots,\alpha_P \in \R$ and $D\in \N\setminus \{0\}$. For all $Q\geq 2$, one can find
 an integer $n \in \{D,\ldots ,DQ^P \}$ such that 
 $$\max_{1\leq j \leq P} \mathrm{dist}(n\alpha_j, \Z) \leq \frac{1}{Q}.$$
\end{lemma}
Applying this lemma with $\alpha_j=\frac{\theta_j}{2\pi}$, for all $\eta>0$ we can find a sequence $N_\ell$ with $N_\ell \rightarrow \infty$ as $\ell$ goes to $\infty$, such that for
all $\ell$, we have
$$\max_j  \left \vert  e^{iN_\ell\theta_j}-1 \right \vert \leq \eta.$$
Therefore we have for all $\ell$, 
$$\vert S(N_\ell) \vert \geq \left \vert \sum_j Q_{j,f,g}(N_\ell) \right \vert-\eta \sum_{j} \vert Q_{j,f,g}(N_\ell) \vert .$$
Now consider the quantity given by
$$R(N):=\frac{ \left \vert \sum_j Q_{j,f,g}(N) \right \vert }{\sum_{j} \vert Q_{j,f,g}(N) \vert}.$$
If we assume that $f,g$ are such that the sum of coefficients
$$ \sum_j \lambda_j^{-d_j +1}L_{m}(g \mathcal{N}_j^{d_j-1} \mathcal{P}_j f)\neq 0,$$
then $R(N)$ makes sense for all $N$ large and $\lim_{N\rightarrow +\infty} R(N)$ exists and is {\it non-vanishing}. We can therefore choose $\eta>0$ such that for all $N$ large we have
$$\eta\leq \frac{1}{2}R(N).$$
The proof is then done because we get for all large $\ell$,
$$\vert C_{\varphi,\psi}(N_\ell)\vert \geq \widetilde{\rho}^{N_\ell} \left ( \frac{1}{2}\vert \sum_j Q_{j,f,g}(N_\ell)\vert +o(1) \right)$$
$$\geq C\widetilde{\rho}^{N_\ell}$$ for some $C>0$. It remains to check that we can adjust $f,g$ such that
$$ T_g(f):=\sum_j \lambda_j^{-d_j +1}L_{m}(g \mathcal{N}_j^{d_j-1} \mathcal{P}_j f)\neq 0.$$
Without loss in generality, we can pick $j_0$ such that $\mathcal{N}_{j_0}^{d_{j_0}-1}\neq 0$ and choose $f \in Im(\mathcal{P}_{j_0})$ with $ \mathcal{N}_{j_0}^{d_{j_0}-1} f \neq 0$. Then we have 
$$T_g(f)=\lambda_{j_0}^{-d_{j_0} +1}L_{m}(g \mathcal{N}_{j_0}^{d_{j_0}-1} f),$$ and by property $3)$ from Proposition \ref{mainprop}, we can choose $g=g_0$ to be a trigonometric polynomial such that $T_{g_0}(f)$ is non vanishing. Consider now the functional $f\mapsto T_{g_0}(f)$. It's now a non trivial continuous linear form on $\mathcal{H}_{r,M}$ and by density
we can choose a trigonometric polynomial $f_0$ such that again $T_{g_0}(f_0)\neq 0$. The proof is done. $\square$ 
\section{Existence of non trivial spectra for $\lt_q$ via frequency averaging} 
In this section, we will prove Proposition \ref{mainprop} and its two statements. The main ideas will revolve around the trace formula
$$\mathrm{Tr}( \lt_q^n)=\sum_{A^nx=x} \frac{e^{2i\pi q \tau^{(n)}(x)}}{\vert \det(I-D_xA^n)\vert},$$
and different ways to estimate (from below and above) this oscillating sum via averaging techniques. We start by a basic a priori bound.
\subsection{An upper bound on the trace}
For all $q$, let $\rho(q,\tau)$ denote the spectral radius of $\lt_q$, i.e.
$$\rho(q,\tau):=\sup_{j}\vert \lambda_j(\lt_q)\vert.$$
\begin{prop}
\label{apriori1}
For all $R\geq 1$, there exists $C(R)$ depending only on $R$ and the map $A$, such that for all $n\in \N$ and $q \in \Z\setminus \{0\}$,
we have
$$\vert \mathrm{Tr}(\lt_q^n)\vert \leq C(R)q^2 \left(\Vert \tau \Vert_{r,\infty}^2+\Vert \tau \Vert_{r,\infty}+1\right)
\max\{e^{-Rn};\rho(q,\tau)^n \}.$$
\end{prop}
\noindent {\it Proof}. We first start to write
$$\vert \mathrm{Tr}(\lt_q^n)\vert \leq \sum_{j=1}^\infty \vert \lambda_j(q)\vert^n$$
$$\leq N (\rho(q))^n+\sum_{j=N+1}^\infty \vert \lambda_j(q)\vert^n,$$
where $N$ will be adjusted later on. Using the a priori estimate for the eigenvalues from Proposition \ref{blackbox}, we have
$$\sum_{j=N+1}^\infty \vert \lambda_j(q)\vert^n\leq e^{n(\alpha+\vert q\vert \Vert \tau \Vert_{r,\infty})} \sum_{j=N+1}^\infty e^{-n \beta \sqrt{j}},$$
where the constant $C$ from Proposition \ref{blackbox} is written as $C=e^{\alpha}$.
On the other hand, we have
$$\sum_{j=N+1}^\infty e^{-n \beta \sqrt{j}}\leq \int_{N}^\infty e^{-n\beta\sqrt{t}}dt=\frac{2}{n^2}e^{-\beta n \sqrt{N}}\left ( \frac{n\sqrt{N}}{\beta} +\frac{1}{\beta^2}  \right).$$
Choosing 
$$N\geq 1+\left( \frac{\alpha+\vert q \vert \Vert \tau \Vert_{r,\infty} +R}{\beta} \right)^2,$$
so that 
$$e^{n(\alpha+\vert q\vert \Vert \tau \Vert_{r,\infty})} \sum_{j=N+1}^\infty e^{-n \beta \sqrt{j}}\leq \frac{2}{n^2}e^{-nR} 
\left ( \frac{n\sqrt{N}}{\beta} +\frac{1}{\beta^2}  \right).$$
We now have obtained
$$\vert \mathrm{Tr}(\lt_q^n)\vert \leq N (\rho(q))^n+\frac{2}{n^2}e^{-nR} 
\left ( \frac{n\sqrt{N}}{\beta} +\frac{1}{\beta^2}  \right),$$
and the proof is done. $\square$
\subsection{Averaging over the frequency $q$}
In this section we shall prove part $1)$ of Proposition \ref{mainprop}. First we need an observation on topological pressure of  the unstable jacobian and weighted sums over periodic orbits that arise from the trace formulas. More precisely we have the following fact.
\begin{lemma}
\label{pressure1}
Let $\sigma>0$.
 For all $\epsilon>0$, one can find a constant $C>0$ such that for all $n$ large enough
 $$C^{-1}e^{n(P(-\sigma\log J^u)-\epsilon)}\leq \sum_{A^n x=x} \frac{1}{\vert \det(I-D_x A^n)\vert^\sigma}\leq C e^{n(P(-\sigma\log J^u)+\epsilon)}.$$
\end{lemma}
\noindent {\it Proof}. We recall that the Anosov structure says that at each point $x\in \T$, we have a splitting 
$$T_x \T=E^u_x \oplus E^s_x,$$
with $D_xA(E_x^u)=E_{Ax}^u$, $D_xA(E_x^s)=E_{Ax}^s$ and there exist constants $C_1,C_2, \lambda^u,\lambda^s >0$ such that for all $n\geq 0$,
$$\Vert D_xA^n\vert_{E_x^s} \Vert\leq C_1 e^{-\lambda^s n},\   \Vert D_xA^{-n}\vert_{E_x^u} \Vert\leq C_2 e^{-\lambda^u n}.$$
Whenever $A^nx=x$, we have two mappings $D_xA^n:E_x^u\rightarrow E_x^u$ and $D_xA^n:E_x^s\rightarrow E_x^s$. Therefore
we have
$$\det(I-D_xA^n)=\det(I-D_xA^n\vert_{E^u_x}) \det(I-D_xA^n\vert_{E^s_x})$$
$$=\det(D_xA^n\vert_{E^u_x})\det(I-D_xA^{-n}\vert_{E^u_x})\det(I-D_xA^n\vert_{E^s_x}).$$
By exponential decay of both $\Vert D_xA^n\vert_{E_x^s} \Vert$ and $\Vert D_xA^{-n}\vert_{E_x^u} \Vert$ as $n\rightarrow +\infty$,
we deduce that there exists a constant $C>0$, such that for all $n$ large, we have
$$ \vert \det(I-D_xA^n)\vert \leq C \vert \det(D_xA^n\vert_{E^u_x}) \vert=CJ^u(x)J^u(Ax)\ldots J^u(A^{n-1}x),$$
so that
$$\sum_{A^n x=x} \frac{1}{\vert \det(I-D_x A^n)\vert^\sigma}\geq C \sum_{A^n x=x} e^{-\sigma(\log J^u)^{(n)}(x)}.$$
It is then a standard fact, that when $A$ is a {\it topologically mixing Anosov map}, we have for all real valued H\"older potential $\varphi$,  $$\lim_{n\rightarrow +\infty} \left ( \sum_{A^n x=x} e^{\varphi^{(n)}(x)} \right )^{1/n}=e^{P(\varphi)}.$$
For references, see the classics \cite{BowenBook,ParryPollicott}. For a more modern treatment, see also \cite{BaladiBook2}, Chapter 7, 
Corollary 7.7. The proof of the lower bound is done. For the upper bound, the exact same ideas work straightforwardly. $\square$

We now proceed toward a proof of proposition \ref{mainprop}, first part. Let us set
$$S(n,q):=\mathrm{Tr}( \lt_q^n)=\sum_{A^nx=x} \frac{e^{2i\pi q \tau^{(n)}(x)}}{\vert \det(I-D_xA^n)\vert}.$$
Our goal is to obtain some decent lower bounds on $\vert S(n,q)\vert$. Pointwise, this is quite a desperate task, but we will rely instead
on an averaged estimate by summing carefully over the frequency parameter $q$. We pick a $C_0^\infty$ test function on $\R$ having the following set of properties:\footnote{The existence of such a test function is a folklore fact. Start with the usual $C_0^\infty$ bump function given by
$\varphi_0(x)=e^{-\frac{1}{1-x^2}}\chi_{[-1,+1]}(x)$.  To make sure that the Fourier transform is
positive consider then the convolution $\psi=\varphi_0 \star \varphi_0$ which obviously has now the desired properties.}
$$\left \{   \begin{array}{c}
\forall x\in \R,\ \psi(x) \geq 0,\mathrm{and }\ \psi(0)>0\\
\mathrm{supp}(\psi)\subset [-2,+2], \mathrm{and}\ \forall \xi \in \R,\ \widehat{\psi}(\xi) \geq 0\\
\end{array}
\right. $$
where $\widehat{\psi}$ is the Fourier transform defined by
$$\widehat{\psi}(\xi):=\int_\R \psi(x)e^{-ix\xi}d\xi.$$ We now set for some $T>0$, 
$$\psi_T(x):=\psi\left(\frac{x}{T}\right),$$
and consider the quantity
$$\sum_{q\in \Z} \psi_T(q)\vert S(n,q)\vert^2= \sum_{A^nx=x\atop A^nx'=x'} \sum_{q \in \Z} \psi_T(q) 
\frac{e^{2i\pi q(\tau^{(n)}(x)-\tau^{(n)}(x'))}}{\vert \det(I-D_xA^n)\vert\vert \det(I-D_{x'}A^n)\vert}.$$
To compute this average, we use {\it Poisson summation formula} 
which says that given a rapidly decaying test function $\varphi \in \mathcal{S}(\R)$,
we have the celebrated identity
$$\sum_{q\in \Z} \varphi(q)=\sum_{k\in \Z} \widehat{\varphi}(2\pi k).$$
Therefore we have 
$$\sum_{q\in \Z} \psi_T(q) e^{2i\pi q(\tau^{(n)}(x)-\tau^{(n)}(x'))}=T\sum_{p\in \Z} 
\widehat{\psi}\left(2\pi T(p-\tau^{(n)}(x)+\tau^{(n)}(x')) \right).$$
Because $\widehat{\psi}$ is {\it positive}, we can obviously bound from below (by forgetting all the non-diagonal terms) for all $p$
$$ \sum_{A^nx=x\atop A^nx'=x'} 
\frac{\widehat{\psi}\left(2\pi T(p-\tau^{(n)}(x)+\tau^{(n)}(x')) \right)}{\vert \det(I-D_xA^n)\vert\vert \det(I-D_{x'}A^n)\vert} \geq \sum_{A^nx=x} 
\frac{\widehat{\psi}\left(2\pi Tp \right)}{\vert \det(I-D_xA^n)\vert^2}.$$
Observe now that because $\widehat{\psi}$ is in the Schwartz class (rapid decay), then we have 
$$\sum_{p\neq 0}\widehat{\psi}(2\pi Tp)=O\left(T^{-1}\right),$$
which tells us that for large $T$, we can as well drop all non zero $p$ terms in the sum $\sum_{p\in \Z}\widehat{\psi}(2\pi Tp)$, so that we end up with the lower bound (we use Lemma \ref{pressure1})
$$\frac{1}{T}\sum_{q\in \Z} \psi_T(q)\vert S(n,q)\vert^2 \geq \sum_{A^nx=x} 
\frac{\widehat{\psi}(0)}{\vert \det(I-D_xA^n)\vert^2}\geq ce^{n(P(-2\log J^u)-\epsilon)},$$
for some $c>0$ and all $n$ large enough. Using Proposition \ref{apriori1}, we have obtained
$$ e^{n(P(-2\log J^u)-\epsilon)}\leq C_{R,\epsilon} \frac{1}{T} \sum_{\vert q\vert\leq 2T}(q^4+1)\max\{e^{-2Rn}, \rho(q)^{2n} \}.$$
identity valid for all $T\geq 1$ and $n$ large. We now fix $0<\eta<1$ and choose $R=\vert \log \eta \vert$.
We end the proof by contradiction. Assume that there exists $q_0>0$ such that for all $\vert q\vert \geq q_0$,
we have $\rho(q)\leq \eta$. We get therefore (using the fact that unconditionally $\rho(q)\leq 1$) 
$$\frac{1}{T} \sum_{\vert q\vert\leq 2T}(q^4+1)\max\{e^{-2Rn}, \rho(q)^{2n} \}
\leq \frac{1}{T} \sum_{\vert q\vert < q_0}(q^4+1)+\frac{\eta^{2n}}{T} \sum_{q_0\leq \vert q\vert\leq 2T}(q^4+1)$$
$$\leq O(T^{-1})+ O(\eta^{2n}T^4).$$
For notational simplicity, we set $P:=\vert P(-2\log J^u) \vert$. We set in the sequel
$$n=[\beta \log T],$$
where $\beta>0$ will be adjusted later on. As $T\rightarrow +\infty$ we have
$$T^{-\beta (P+\epsilon)}=O(T^{-1})+O(T^{2\beta\log \eta +4}).$$
we get a contradiction whenever
$$\beta P<1\ \mathrm{and}\ \beta P<2\beta \vert \log \eta \vert -4,$$
which leads to choose 
$$\eta=e^{\frac{5}{2}P(-2\log J^u)-\overline{\epsilon}},$$
for all $\overline{\epsilon}>0$, and the proof is done. 
\section{An improved lower bound via a probabilistic technique}
The goal of this section is to show how to improve the lower bound of Theorem \ref{main}, 1) via a different argument. Instead of averaging over the frequency parameter $q$, we will consider some random  "coupling functions" $\tau$ of the form
$$\tau(z)=P_{N}(z),$$
where $P_{N}(z)$ belongs to a suitable ensemble of {\it random trigonometric polynomials $\mathcal{P}_N$}. 
The game is to estimate from below the expectation 
$\E(\vert \mathrm{Tr}(\lt_{q,\tau}^n)\vert^2)$, where $\mathrm{Tr}(\lt_{q,\tau}^n)$ is seen as a random variable. This technique will allow us to overcome the "exponent loss", artefact of the frequency averaging technique, and will rely also on a positivity argument. 
We will also prove that rapid mixing occurs with probability $1$ in $\mathcal{P}_N$, which will occupy section $\S 4.2$.

\subsection{The set of random trigonometric polynomials $\mathcal{P}_N$}
 Our goal is to define a set of random {\it real valued} trigonometric polynomials. We could use a Fourier basis of real valued trigonometric functions,
based on product of sines and cosines, but that would lead to cumbersome notations. Instead, we choose a more conceptual route using a fixed basis of {\it real valued} eigenfunctions of the Laplacian.
Let $\Delta:=-\partial_1^2-\partial_2^2$ be the flat Laplacian on $\T$.
We choose an $L^2(\T)$ basis of {\it real eigenfunctions} $\varphi_0,\varphi_1,\ldots,\varphi_j$ of the Laplacian such that
$$\Delta \varphi_j=\lambda_j \varphi_j,$$
where the eigenvalues $\lambda_0=0<\lambda_1\leq \lambda_2\leq \ldots \leq \lambda_j$ are repeated according to multiplicity.
Each eigenfunction $\varphi_j(z)$ is a trigonometric polynomial of the form
$$\varphi_j(z)=\sum_{\alpha \in \Z^2\atop \lambda_j=4\pi^2\Vert \alpha \Vert^2} C_\alpha e^{2i\pi \alpha.z},$$
with $\Vert \alpha \Vert^2=\alpha_1^2+\alpha_2^2$, and $C_{-\alpha}=\overline{C_\alpha}$ are coefficients subject to the $L^2$ normalization. We will need to use two basic facts on these eigenfunctions:
\begin{itemize}
 \item (Weyl law). As $N\rightarrow +\infty$, we have $\# \{ j\ :\ \vert \lambda_j \vert \leq N \}=O(N^2)$.
 \item ($L^\infty$-growth). As $j\rightarrow \infty$, we have $\Vert \varphi_j \Vert_\infty=O_\epsilon(\lambda_j^\epsilon)$, for all $\epsilon>0$.
\end{itemize}
The first fact is easy and follows from a crude upper bound on the number of lattice points in a disc. The other claim
is less trivial but standard and follows from estimating the number of lattice points on a circle, which in turn is related to the number of representations of a given integer as a sum of two squares.

Let $(\Omega,\Prob)$ be a probability space,
and assume that there exists a sequence $(X_j)_{j\in \N}$ of real valued, independent random variables on $\Omega$, whose common probability law is Gaussian centered with variance $1$, i.e. each $X_j$ obeys the law (for all measurable $A$)
$$\Prob(X_j \in A)=\frac{1}{\sqrt{2\pi}}\int_A e^{-\frac{x^2}{2}}dx.$$
Fix $N$ a large integer and consider the random trigonometric polynomial given by
$$P_N(z)=\sum_{j\ :\ \lambda_j \leq N} X_j \varphi_j(z),$$
and this random ensemble is denoted by $\mathcal{P}_N$.
We start by a basic Lemma.
\begin{lemma}
\label{moments} 
For all $r>0$, for all polynomial $Q(x)=\sum_j a_j x^j$, with $a_j \geq 0$, the expectation 
$$\E\left (  Q(\Vert P_N \Vert_{r,\infty})  \right)<+\infty.$$
\end{lemma}
\noindent {\it Proof}. Let us prove that the random variable $\Vert P_N \Vert_{r,\infty}$ has finite even moments to all orders. 
By Schwarz inequality, we have ($x\in \T$, $y\in [-r,+r]^2$),
$$\vert P_N(x+iy)\vert^2\leq K_{N,r} \sum_{\lambda_j\leq N}{X_j^2},$$
where we have set 
$$K_{N,r}:=\sup_{x\in \T,y\in [-r,+r]^2} \left( \sum_{\lambda_j\leq N} \vert \varphi_j(x+iy)\vert^2 \right).$$
Consequently, we have 
$$\vert P_N(x+iy)\vert^{2p}\leq K_{N,r}^p \sum_{j_1,j_2,\ldots,j_p} X_{j_1}^2 X_{j_2}^2 \ldots X_{j_p}^2.$$
Because the random gaussian variables $X_j$ are assumed to be independent and have finite moments to all orders, we can conclude 
(without having to compute it) that
$$\sum_{j_1,j_2,\ldots,j_p} \E(X_{j_1}^2 X_{j_2}^2 \ldots X_{j_p}^2)<+\infty.$$
Using Schwarz inequality, we deduce now that all the odds moments are also finite, and the proof is done. $\square$

\subsection{Rapid mixing occurs almost surely in $\mathcal{P}_N$.} In the following, we prove that rapid mixing is almost sure in $\mathcal{P}_N$. This part
of the paper will definitely not surprise the experts, but there are nevertheless some technical details that have to be addressed.

Here we recall that the real analytic Anosov map $A:\T\rightarrow\T$ is
assumed to be topologically mixing, which definitely occurs if $A$ is $C^1$-close to a linear hyperbolic map. The proof is twofold: first we use symbolic dynamics
to show that rapid mixing follows from an estimate on transfer operators due to Dolgopyat \cite{Dolgopyat}. This estimate holds under a diophantine
condition satisfied by the "roof" function $\tau$. We then show that when $\tau \in \mathcal{P}_N$, this diophantine hypothesis holds with full probability.

According to Bowen \cite{BowenBook}, using Markov partitions, there exists a topologically mixing subshift of finite type $(\Sigma,\sigma)$, and a continuous map
$$\Pi: \Sigma \rightarrow \T,$$
such that $A\circ \Pi=\Pi\circ \sigma$. We recall that 
$$\Sigma=\{ (x_i)\in \{1,\ldots,p\}^\Z\ :\ \forall i\in \Z,\ M(x_i,x_{i+1})=1 \},$$
where $M$ is some $p\times p$ aperiodic matrix whose entries are $0$ or $1$, while
the shift map $\sigma$ is given by
$$\sigma(\xi)_j=\xi_{j+1}.$$ 
This compact product space $\Sigma$ can be equiped with an ultrametric distance $d_\theta$, with
$0<\theta<1$ given by
$$d_\theta(x,\xi):=\left \{ 1\ \mathrm{if}\ x_0\neq \xi_0 \atop \theta^ {1+\max\{k\geq 0\  :\ x_j=\xi_j\ \forall\ \vert j \vert \leq k \}}. \right. $$
For an adequate choice of $\theta$, the map $\Pi$ becomes Lipschitz continuous. Similarly, the one-sided subshift $\Sigma^+$ is simply
$$\Sigma^+=\{ (x_i)\in \{1,\ldots,p\}^\N\ :\ \forall i\in \N,\ M(x_i,x_{i+1})=1 \},$$
endowed with the same metric $d_\theta$.
We denote by $\ft$ and $\ftp$, the function spaces of Lipschitz continuous functions on $\Sigma,\Sigma^+$ endowed with 
$$\Vert f\Vert_{\ft}:=\Vert f\Vert_\infty +\sup_{x\neq \xi} \frac{\vert f(x)-f(\xi)\vert}{d_\theta (x,\xi)}.$$ 
The SRB measure $\mu_{SRB}$ on $\T$ is the pull back of the $\sigma$-invariant equilibrium measure $\widetilde{\mu}$ on $\Sigma$
associated to the Holder potential $-\log J^u$. On the one-sided subshift $\Sigma^+$ there is an associated measure (called again $\widetilde{\mu}$)
wich satisfies
$$\int_\Sigma f d\widetilde{\mu}= \int_{\Sigma^+} f d\widetilde{\mu},$$
whenever $f \in C^0(\Sigma)$ depends only on positive coordinates (so that $f$ corresponds to an element of $C^0(\Sigma^+)$).
The map $\sigma:\Sigma^+\rightarrow \Sigma^+$ acts on $L^2(\Sigma^+,d\widetilde{\mu})$ and there is a unique adjoint transfer operator $L_h$
such that 
$$\int f\circ \sigma g d\widetilde{\mu}=\int f L_h (g) d\widetilde{\mu}$$
given by
$$L_h(g)(x)=\sum_{\sigma y= x} e^{h(y)}g(y),$$
for some $h\in \ftp$. The estimate needed to prove rapid decay is from Dolgopyat \cite{Dolgopyat}, see also \cite{Naud2} for the version stated below.
\begin{prop}
 \label{Dolg}
 Let $\tau \in \ftp$. Assume that for some $n\geq 1$, there exist two periodic points $\sigma^n\xi=\xi$, $\sigma^n x=x$ such that the vector
 $$(\tau^{n}(x),\tau^{n}(\xi))\in \R^2$$
is diophantine. Then there exists $C>0,\beta>0,\gamma>0$ such that for all $q\neq 0$, the transfer operator
 $$\lp_q(f):=\sum_{\sigma y= x} e^{h(y)+2i\pi q \tau(y)}f(y)$$
 satisfies on $\ftp$ for all $q\neq 0$ and $N\geq 0$, 
 $$\Vert \lp_q^N\Vert_{\ftp}\leq C \vert q\vert^\beta \left( 1-\frac{1}{\vert q\vert^\gamma} \right)^N.$$
\end{prop}
Let us now show briefly how this spectral estimate implies rapid mixing. We start with two smooth observables $F,G\in C^\infty(\T)$, The correlation function then writes as
$$C_{F,G}(N)=\int_{S^1}\int_{\Sigma} F(A^Nx,\tau^{(N)}(x)+\omega)G(x,\omega)d\widetilde{\mu}(x)d\omega.$$
where all the functions involved are pulled back to $\Sigma$ via $\Pi$. The next step is to remark that one can assume that $\tau$ depends only on future coordinates via a well known lemma,
see for example \cite{BaladiBook1}, lemma 1.3. Indeed there exists $\varphi \in \mathcal{F}_{\sqrt{\theta}}$ and $\tau^+\in \mathcal{F}^+_{\sqrt{\theta}}$ such that
$$\tau=\tau^++\varphi-\varphi\circ \sigma.$$
A change of variable then yields
$$ C_{F,G}(N)=\int_{S^1}\int_{\Sigma} \widetilde{F}(A^Nx,(\tau^+)^{(N)}(x)+\omega)\widetilde{G}(x,\omega)d\widetilde{\mu}(x)d\omega,$$
where $\widetilde{F}(x,\omega)=F(x,-\varphi(x)+\omega)$, $\widetilde{G}(x,\omega)=G(x,-\varphi(x)+\omega)$.
To avoid more complicated notations, we will still use the notation $\ftp$, even though it is clear from the precedent remark
that a loss of regularity is required to be able to reduce the problem to the one-sided subshift $\Sigma^+$.
A more delicate argument (see Dolgopyat \cite{Dolgopyat}, section 2.3) involving projections on "unstable manifolds" allows to approximate simultaneously $F,G$ at exponential rate by functions that depend 
only on future coordinates thus reducing the problem to the case of $\Sigma^+$. Assume from now on that $\tau\in \ftp$ and that for all $\omega \in S^1$, 
$F(.,\omega), G(.,\omega)\in \ftp$. Assume that the $\omega$ dependence is $C^\infty$, and set for all $q\in \Z$,
$$\widehat{F}_q(x):=\int_{S^1} F(x,\omega)e^{-2i\pi q \omega}d\omega.$$
By repeated integration by parts, we have for all $k\geq 0$, $q\neq 0$, 
$$\Vert \widehat{F}_q \Vert_{\ftp}\leq B_k \vert q \vert^{-k}.$$
The Correlation functions now expands as
$$C_{F,G}(N)=\sum_{q \in \Z} \int_{\Sigma^+} (\widehat{F}_q\circ \sigma^N) e^{2i\pi q (\tau^+)^{(N)}}\widehat{G}_{-q}d\widetilde{\mu}$$
$$=\int_{\Sigma^+} (\widehat{F}_0\circ \sigma^N) \widehat{G}_{0}d\widetilde{\mu}+\sum_{q \neq 0} \int_{\Sigma^+} \widehat{F}_q  \lp_q^N(\widehat{G}_{-q})d\widetilde{\mu}$$
$$=O(\eta^N)+O_k \left(\sum_{q\neq 0} \vert q \vert^{-2k+\beta} e^{-\frac{N}{\vert q \vert^\gamma}} \right). $$
On the other hand,
$$\sum_{q\neq 0} \vert q \vert^{-2k+\beta} e^{-\frac{N}{\vert q \vert^\gamma}}=
O( e^{-\frac{N}{\vert \widetilde{N} \vert^\gamma}})+O\left(\sum_{\vert q\vert \geq \widetilde{N}} \frac{1}{\vert q\vert^{2k-\beta}}\right).$$ Choosing $\widetilde{N}=[ N^{1/\gamma-\epsilon}]$ proves rapid mixing since $k$ can be taken arbitrarily
large.

We have now to show that whenever $\tau\in \mathcal{P}_N$, the diophantine condition from Proposition \ref{Dolg} is satisfied almost surely. Remark that this condition is on $\tau^+$, but because of the cohomological relation it is enough to prove it for the original $\tau$. We recall that a vector $(x_1,x_2)\in \R^2$ is diophantine if and only if there exists $m_0>0$ such that for all $\alpha \in \Z^2$ with $\vert \alpha \vert=\vert \alpha_1\vert +\vert \alpha_2 \vert\geq 2$, we have
$$\vert \alpha_1 x_1+\alpha_2 x_2 \vert \geq \vert \alpha \vert^{-m_0}.$$
Let us pick $n$ large enough and two $n$-periodic points $x,y\in \Sigma$ and denote, slightly abusing notations,
their projections $\Pi(x),\Pi(y) \in \T$ again by $x,y$. We assume that the orbits
$$\{ x,Ax,\ldots,A^{n-1}x\},\  \{ x,Ax,\ldots,A^{n-1}x\}$$
are different, which is of course possible by taking $n$ large enough. Given $\alpha\in \Z^2$, we will set 
$$Y_\alpha:=\alpha_1P_N(x)^{(n)}+\alpha_2 P_N(y)^{(n)}=
\sum_{\lambda_j\leq N} X_j \left ( \alpha_1\varphi_j(x)^{(n)}+\alpha_2 \varphi_j(y)^{(n)}\right).$$
Then $Y_\alpha$ is a Gaussian variable with expectation $0$ and variance
$$\sigma^2(\alpha)=\sum_{\lambda_j\leq N} \left ( \alpha_1\varphi_j(x)^{(n)}+\alpha_2 \varphi_j(y)^{(n)}\right)^2.$$
We need a lower bound on the variance $\sigma^2(\alpha)$ which is given by the following Lemma.
\begin{lemma}
\label{deviation}
Using the above notations, there exists $C>0$ and $N_0$ such that for all $N\geq N_0$ and all $\alpha\in \Z^2$ with $\vert \alpha\vert\geq 2$, we have
$$\sigma^2(\alpha)\geq C\vert \alpha\vert^{2}.$$
\end{lemma}
\noindent {\it Proof}. Because both periodic orbits $\{ x,Ax,\ldots,A^{n-1}x\}$ and $\{ y,Ay,\ldots,A^{n-1}y\}$ are different,
there exists a $C^\infty$ real valued function $\psi$ on $\T$ such that $\psi(x)=\psi(Ax)=\ldots=\psi(A^{n-1}x)=\mathrm{sign}(\alpha_1)$ and $\psi(y)=\psi(Ay)=\ldots=\psi(A^{n-1}y)=\mathrm{sign}(\alpha_2)$, where $\mathrm{sign}(x)=\frac{x}{\vert x \vert}$ if $x\neq 0$. Therefore we have
$$ \alpha_1\psi(x)^{(n)}+\alpha_2 \psi(y)^{(n)}=n\vert \alpha \vert.$$
On the other hand, being smooth, $\psi$ has a uniformly convergent expansion in the eigenfunction basis $(\varphi_j)$ as 
$$\psi(z)=\sum_{j=0}^\infty C_j \varphi_j(z),$$
where $C_j$ is given by
$$C_j=\int_{\T} \psi(x)\varphi_j(x)dm(x).$$
Applying several times the flat Laplacian $\delta$ and using Green's formula shows that the coefficients $C_j$ are rapidly
decreasing i.e. for all $k\geq 0$, there exists $B_k$ such that 
$$\vert C_j\vert\leq B_k \lambda_j^{-k}.$$
Therefore we can write (fixing $k$ large enough and $\epsilon>0$ small)
$$n\vert \alpha\vert=\vert \alpha_1\psi(x)^{(n)}+\alpha_2 \psi(y)^{(n)}\vert\leq \sum_{\lambda_j \leq N} \vert C_j\vert  
\vert \alpha_1\varphi_j(x)^{(n)}+\alpha_2 \varphi_j(y)^{(n)}\vert+
M_{k,\epsilon}\vert \alpha\vert \sum_{\lambda_j >N} \lambda_j^{-k+\epsilon}.$$
It is now clear that by taking $N$ large enough, we have
$$ \frac{1}{2} \vert \alpha \vert \leq \sum_{\lambda_j \leq N} \vert C_j\vert  
\vert \alpha_1\varphi_j(x)^{(n)}+\alpha_2 \varphi_j(y)^{(n)}\vert,$$
and the proof is done by applying Schwarz inequality. $\square$.

We now estimate the probability that the random vector $(\tau^{(n)}(x),\tau^{(n)}(y))$ is {\it not} diophantine, $x,y,n$ being fixed. According to our definition, $(\tau^{(n)}(x),\tau^{(n)}(y))$ being not diophantine corresponds to the event
$$\bigcap_{m_0\geq 1} \bigcup_{\vert \alpha \vert \geq 2} 
\left \{ \vert Y_\alpha \vert < \vert \alpha \vert^{-m_0}\right \}.$$
Since we have
$$\Prob( \vert Y_\alpha \vert < \vert \alpha \vert^{-m_0})\leq 
\frac{1}{\sigma(\alpha)\sqrt{2\pi}} \int_{-\vert \alpha \vert^{-m_0}}^{+\vert \alpha \vert^{-m_0}}
e^{-\frac{x^2}{2\sigma^2(\alpha)}}dx\leq \frac{2\vert \alpha \vert^{-m_0}}{\sigma(\alpha)\sqrt{2\pi}},$$
we can use Lemma \ref{deviation} to write
$$ \Prob(\bigcup_{\vert \alpha \vert \geq 2} 
\left \{ \vert Y_\alpha \vert < \vert \alpha \vert^{-m_0}\right \})\leq 
C\sum_{\vert \alpha \vert \leq 2} \vert \alpha \vert^{-m_0-1}.$$
On the other hand, for all $m_0\geq 2$, we have
$$\sum_{\alpha \in \Z^2, \vert \alpha \vert\geq 2} \vert \alpha \vert^{-m_0-1}=O(2^{-m_0+1}),$$
and therefore,
$$\Prob(\bigcap_{m_0\geq 1} \bigcup_{\vert \alpha \vert \geq 2} 
\left \{ \vert Y_\alpha \vert < \vert \alpha \vert^{-m_0}\right \})=\lim_{m_0\rightarrow +\infty} 
\Prob(\bigcup_{\vert \alpha \vert \geq 2} 
\left \{ \vert Y_\alpha \vert < \vert \alpha \vert^{-m_0}\right \})=0.$$
The claim on almost sure rapid mixing is proved. We now conclude this section by explaining why almost sure exponential mixing
is likely to be much harder. Because we are in a low dimensional hyperbolic case, one can use results from Palis-Takens \cite{PalisTakens}, which show that the stable/unstable foliation for the anosov map $A$ is of class $C^{1+\alpha}$, for some $\alpha>0$.
One can then "quotient out" the stable direction to reduce the problem to a purely expanding situation and try to use Dolgopyat's estimates in a $C^{1+\alpha}$ setting, as in the recent work of Butterley-War \cite{BW1} for Anosov flows. There is however a major issue
there: unlike in the flow case, our roof function $\tau$ does not come naturally as a return time function and is not constant on stable leaves. Projecting $\tau$ on unstable leaves would produce a merely H\"older function, destroying all options to use \cite{BW1}.
On the other hand, Liverani as shown in \cite{Liverani} that by working directly with Anisotropic norms, one can avoid this situation.
There is however a major ingredient in \cite{Liverani} that we cannot use here: it's the contact structure that confers to the "return time"
some very strong non integrability properties. As a conclusion, it seems to us that Tsujii's recent work \cite{Tsujii} might be the way
to go to prove generic exponential mixing in this setting.

\subsection{Spectral lower bound via expectation}
We can now give a proof of Proposition \ref{mainprop}. We fix $N$ large enough as in the previous section such that with full probability,
rapid mixing holds, otherwise our result would be trivial. Indeed, any roof function $\tau$ which is cohomologous to a constant function
produces a non-mixing extension for which the claimed lower bound trivially holds.
We fix some $\varepsilon>0$. By the trace formula, and using independence of the gaussian variables $X_j$, we have 
$$\E(\vert \mathrm{Tr}(\lt_{q,P_N})\vert^2)=\sum_{A^nx=x\atop A^nx'=x'} 
\frac{\E \left(e^{2i\pi q (P_N^{(n)}(x)-P_N^{(n)}(x'))}\right)}{\vert \det(I-D_xA^n)\vert\vert \det(I-D_{x'}A^n)\vert}$$
$$=\sum_{A^nx=x\atop A^nx'=x'} 
\frac{\prod_{\lambda_j\leq N}\E \left(e^{2i\pi q (\varphi_j^{(n)}(x)-\varphi_j^{(n)}(x'))}\right)}{\vert \det(I-D_xA^n)\vert\vert \det(I-D_{x'}A^n)\vert}.$$   
The characteristic function $\E(e^{-i\xi X})$ of a normal variable (with mean $0$ and deviation $1$) is given by the Fourier transform of the normal distribution which is 
$$ \E(e^{-i\xi X})=e^{-\frac{\xi^2}{2}},$$
and leads to the formula
$$\E(\vert \mathrm{Tr}(\lt_{q,P_N})\vert^2)=\sum_{A^nx=x\atop A^nx'=x'} 
\frac{\exp\left (-2\pi^2q^2\sigma^2_n(x,x')\right)}{\vert \det(I-D_xA^n)\vert\vert \det(I-D_{x'}A^n)\vert},$$
where
$$\sigma^2_n(x,x')=\sum_{\lambda_j\leq N} \left (\varphi_j^{(n)}(x)-\varphi_j^{(n)}(x')\right)^2 .$$
Since this is a sum of positive terms, we can drop all the non diagonal terms and use Lemma \ref{pressure1} to write for all $n$ large,
$$\E(\vert \mathrm{Tr}(\lt_{q,P_N})\vert^2)\geq Ce^{n(P(-2\log J^u)-\varepsilon)}.$$ Fix some $0<\rho_0<1$ and assume now that for all $q\geq q_0$, we have 
$$\Prob( \rho(\lt_{q,P_N})\leq \rho_0)=1.$$
Using Lemma \ref{apriori1} with $R=\rho_0$, we have for all $n$, 
$$\E(\vert \mathrm{Tr}(\lt_{q,P_N}^n)\vert^2)\leq C q^4 \rho_0^{2n}\E\left(  (\Vert \tau \Vert_{r,\infty}^2+\Vert \tau \Vert_{r,\infty}+1)^2 \right).$$
Thanks to Lemma \ref{moments}, we now that
$$\E\left(  (\Vert \tau \Vert_{r,\infty}^2+\Vert \tau \Vert_{r,\infty}+1)^2 \right)<+\infty.$$ 
We have therefore obtained that for all $q\geq q_0$ and $n$ large, 
$$e^{n(P(-2\log J^u)-\varepsilon)}\leq Mq^4 \rho_0^{2n},$$
For some constant $M>0$. We now fix $q\geq q_0$. Since
$$(Mq^4)^{\frac{1}{2n}}\rightarrow 1$$
as $n\rightarrow +\infty$, we choose $n$ large enough so that
$$ (Mq^4)^{\frac{1}{2n}}\leq e^{\frac{\varepsilon}{2}}.$$
By elevating both sides of the inequality to the power $\frac{1}{2n}$, we obtain
$$e^{\frac{1}{2}P(-2\log J^u)-\varepsilon}\leq \rho_0,$$
a contradiction if we take 
$$\rho_0=e^{\frac{1}{2}P(-2\log J^u)-2\varepsilon}.$$ As a conclusion, there exists $q\geq q_0$ such that 
$$\Prob\left( \rho(\lt_{q,P_N})> e^{\frac{1}{2}P(-2\log J^u)-2\varepsilon}\right)>0,$$
which ends the proof of the last claim in the main theorem.

\section{The Hilbert space $\Aniso$ and a priori bounds on the eigenvalues of $\lt_q$ }
In this section we provide the definitions and proofs of Proposition \ref{blackbox} about the main function space $\Aniso$ and
the spectral properties of $\lt_q$ acting on this space. 
\subsection{The anisotropic space $\Aniso$} We use the notation $e_\alpha(z)$, $\alpha \in \Z^2$ for the usual Fourier basis of $L^2(\T)$,
$$e_\alpha(z)=e^{2i\pi\alpha.z}.$$
We will mostly work on the universal cover $\R^2$ of $\T$ and identify functions on $\T$ as $\Z^2$-periodic functions on $\R^2$.
Given $r>0$, we denote by $H_r^2(\T)$ the space of functions $f$ on $\R^2$ which are $\Z^2$-periodic and enjoy a holomorphic extension to
$\R^2+i(-r,+r)^2 \subset \C^2$, and such that the Hardy norm
$$\Vert f\Vert^2_{H^2_r}:=\sup_{y\in (-r,+r)^2} \int_{\T}\vert f(x+iy)\vert^2dm(x)$$
is finite. One can view this space as the Hardy space of holomorphic functions on a Grauert tube of radius $r$ "around" the torus $\T$. 
This is a Hilbert space, and 
$$\left \{ e_\alpha(z)e^{-2\pi r\vert \alpha \vert}  \right \}_{\alpha \in \Z^2},$$
is an orthonormal basis of $H^2_r(\T)$. If we denote by $\widehat{f}(\alpha)$ the Fourier coefficients given by
$$\widehat{f}(\alpha):=\int_{\T}f(x)e_{-\alpha}(x)dm(x),$$
then for all $f\in H^2_r(\T)$, we have the plancherel formula \footnote{The proof follows by computing the inner product 
$\langle f, e_\alpha e^{-2\pi r \vert \alpha \vert} \rangle_{H^2_r}$ and using contour deformation.}
$$\Vert f \Vert^2_{H^2_r}=\sum_{\alpha \in \Z^2 } \vert\widehat{f}(\alpha)\vert^2 e^{4\pi r \vert \alpha \vert}.$$
This Plancherel formula shows that one can also think of $H^2_r(\T)$ as a Sobolev space on $\T$, with exponential weights.

\bigskip
Let $M\in SL_2(\Z)$ be a hyperbolic matrix whose eigenvalues $\mu_M,\mu_M^{-1}$ do satisfy 
$$\vert \mu_M \vert <1,\ \mathrm{while}\ \vert \mu_M^{-1}\vert <1.$$
Let $Id_{\R^2}=P^++P^-$ where $P^\pm$ are the linear projectors on the eigenspaces
$$\mathrm{Ker}( M-\mu^{\pm}Id).$$
Slightly abusing notations, we will write the decomposition for all $x\in \R^2$,
$$x=P^+(x)+P^-(x)=x^++x^-.$$
Notice that for all $x=(x_1,x_2)\in \R^2$, we have by triangle inequality
$$\vert x \vert:=\vert x_1\vert +\vert x_2\vert\leq \vert x^+\vert+\vert x^-\vert,$$
while
$$\vert x^+\vert+\vert x^-\vert \leq C(M) \vert x \vert,$$
where $C(M)>0$ depends only on $M$.
Given a trigonometric polynomial
$$f(z)=\sum_\alpha a_\alpha e_\alpha(z),$$
we set 
$$A_{r,M}(f)(z):=\sum_{\alpha} a_\alpha e^{-2\pi r (\vert \alpha^+\vert-\vert \alpha^-\vert)} e_\alpha(z).$$
As proved in \cite{AlexAdam}, the linear map $A_{r,M}$ has a continuation
$$A_{r,M}:H_r^2(\T)\rightarrow L^2(\T),$$
which is a {\it bounded} linear operator. The function space $\Aniso$ is then defined as the {\it completion} of 
$H^2_r(\T)$ for the norm
$$\Vert f \Vert_{r,M}:=\Vert A_{r,M} f \Vert_{L^2}.$$
A basis of $\Aniso$ is then given by $\{ \rho_\alpha(z) \}_{\alpha \in \Z^2}$ defined by
$$\rho_\alpha(z):=e^{2\pi r(\vert \alpha^+\vert-\vert \alpha^-\vert)} e_\alpha(z).$$
In particular, elements of $\Aniso$ can be formally identified with combinations
$$\sum_{\alpha \in \Z^2} a_\alpha \rho_\alpha(z),$$
with $\sum_\alpha \vert a_\alpha \vert^2<+\infty$. Beware that $\Aniso$ contains distributions with infinite order. Just to give an idea of how wild
some elements of $\Aniso$ can be, let us write $\alpha^-=P^-(\alpha)=(\alpha.\beta^-)\gamma^-$ where $\gamma^-$ is chosen so that $\vert \gamma^-\vert=1$ and $\beta^-$ is some
vector. Then one can check that the formal sum
$$T_{\beta^-}:=\sum_{\alpha\in \Z^2} e^{-\vert \alpha^-\vert^2}e_\alpha=\sum_{\alpha\in \Z^2} e^{- (\alpha.\beta^-)^2}e_\alpha,$$
which belongs to $\Aniso$, is equal formally to the distribution (with infinite order)
$$\sum_{n=0}^\infty \frac{1}{(2\pi)^n n!} \partial_{\beta^-}^{2n}  \delta_0,$$
where $\delta_0$ is the dirac mass at $0$ in $\T$. Notice that $\delta_0$ itself does not belong to $\Aniso$.

\subsection{Convolution operators on $\Aniso$}
In this section, we study the boundedness of multiplication by a given analytic function on $\Aniso$.  More precisely, we will prove the following statement which is enough for the applications we have in mind.
\begin{prop}
 \label{conv}
 There exists $K(M)\geq1$ such that for all $\widetilde{r}>K(M)r$, for all $F \in H^2_{\widetilde{r}}(\T)$, the multiplication operator
 $$T_F:\left \{ \Aniso \rightarrow \Aniso \atop \varphi\mapsto \varphi \times F  \right.$$
 is well defined and bounded. Moreover, there exists $L(r,\widetilde{r},M)>0$ such that 
 $$\Vert T_F \Vert \leq L \Vert F \Vert_{H^2_{\widetilde{r}}}.$$ 
\end{prop}
\noindent {\it Proof}. Let us fix $\widetilde{r}>r>0$. Let $\varphi \in H^2_r$ and $F \in H^2_{\widetilde{r}}$. Obviously, since $F$ has a holomorphic continuation to a domain $\R^2+i(-\widetilde{r},+\widetilde{r})^2$ with $\widetilde{r}>r$, the map 
$$\left \{ H^2_r \rightarrow H^2_r \atop \varphi\mapsto \varphi \times F  \right. $$
is well defined. We need to estimate 
$$\Vert \varphi F \Vert_{r,M}^2=\sum_{\alpha \in \Z^2} \vert \widehat{F\varphi}(\alpha)\vert^2 e^{-4\pi r(\vert \alpha^+\vert -\vert \alpha^-\vert)}.$$
At the Fourier level, multiplication becomes a convolution so that
$$ \widehat{F\varphi}(\alpha)=\sum_{\beta \in \Z^2}\widehat{F}(\beta)\widehat{\varphi}(\alpha-\beta).$$
We can use Schwarz inequality to write
$$  \vert \widehat{F\varphi}(\alpha)\vert^2\leq \left ( \sum_\beta \vert \widehat{F}(\beta) \vert^2 e^{4\pi \widetilde{r} \vert \beta \vert}   \right)
\left ( \sum_\beta  e^{-4\pi \widetilde{r} \vert \beta \vert}  \vert \widehat{\varphi}(\alpha-\beta) \vert^2 \right)$$
$$\leq \Vert F \Vert^2_{H^2_{\widetilde{r}}} \left ( \sum_\beta  e^{-4\pi \widetilde{r} \vert \beta \vert}  \vert \widehat{\varphi}(\alpha-\beta) \vert^2 \right).$$
Using Fubini and a change of variable, we have
$$\sum_\alpha \sum_\beta  e^{-4\pi r ( \vert \alpha^+\vert-\vert \alpha^-\vert)}e^{-4\pi \widetilde{r} \vert \beta \vert}  \vert \widehat{\varphi}(\alpha-\beta) \vert^2=\sum_\gamma \sum_\alpha e^{-4\pi r ( \vert \alpha^+\vert-\vert \alpha^-\vert)}e^{-4\pi \widetilde{r} \vert \alpha-\gamma\vert}  \vert \widehat{\varphi}(\gamma) \vert^2.$$
Since we have
$$\vert \alpha-\gamma\vert \geq C(M)^{-1} (\vert \alpha^+-\gamma^+\vert+\vert \alpha^- -\gamma^-\vert),$$
we can choose $\widetilde{r}$ such that $C(M)^{-1}\widetilde{r}=r+\epsilon$ for some $\epsilon>0$. We hence get by triangle inequality
$$e^{-4\pi r ( \vert \alpha^+\vert-\vert \alpha^-\vert)}e^{-4\pi \widetilde{r} \vert \alpha-\gamma\vert}\leq e^{-4\pi r \vert \gamma^+\vert-\vert \gamma^-\vert)}e^{-4\pi \epsilon \vert \alpha-\gamma \vert}.$$
We have reached
$$ \Vert \varphi F \Vert_{r,M}^2\leq \left( \sum_\alpha e^{-4\pi\epsilon \vert \alpha \vert} \right) \Vert F \Vert^2_{H^2_{\widetilde{r}}}
\Vert \varphi \Vert_{r,M}^2,$$
and the proof is done. $\square$
\subsection{Singular and eigenvalue estimates for $\lt_q$} 
The main result from \cite{AlexAdam} that we need is the following.
\begin{theorem} 
\label{Adam}
Let $M\in SL_2(\Z)$ be a hyperbolic map as above. Assume that $A$ is a real analytic Anosov map with is $C^1$- close enough to $M$. Then for all $r>0$ small enough, the transfer operator (or Koopman operator) $\lt_0: \Aniso\rightarrow \Aniso$ given by
$$\lt_0(\varphi):=\varphi\circ A ,$$
is a well defined compact operator. Moreover there exists $r_0>0$ (depending on $r$) such that for all $\alpha,\beta \in \Z^2$, we have
$$\langle \lt_0( \rho_\alpha),\rho_\beta \rangle_{\Aniso} \leq e^{-r_0(\vert \alpha\vert +\vert\beta\vert)}.$$
\end{theorem}
This theorem implies in particular that $\lt_0$ is a trace class operator. In the following, we fix $\tau:\T\rightarrow \R$ to be a real analytic function
on the torus, and we assume that $r>0$ is fixed small enough so that we can apply Proposition \ref{conv}. Under these assumptions, the transfer operators $\lt_q=e^{2i\pi q \tau} \lt_0$ are all trace class operators on $\Aniso$, which is claim $(1)$ of Theorem \ref{blackbox}. We now prove
claim $(3)$. 

First we need to recall some basic facts about singular values. Our basic reference for that matter is the book \cite{Bsimon}. If $T:\mathcal{H}\rightarrow \mathcal{H}$ is a compact operator acting on a Hilbert
space $\mathcal{H}$, the {\it singular value sequence} is by definition the sequence $\mu_1(T)=\Vert T\Vert\geq \mu_2(T)\geq\ldots \geq \mu_n(T)$
of the eigenvalues of the positive self-adjoint operator $\sqrt{T^*T}$. Our main tool to estimate singular values is provided by the following Lemma.
\begin{lemma}
\label{singular}
Assume that $(e_j)_{j\in J}$ is a Hilbert basis of $\mathcal{H}$, indexed by a countable set $J$. Let $T$ be a compact operator on $\mathcal{H}$.
Then for all subset $I\subset J$ with $\# I=n$ we have 
$$\mu_{n+1}(T)\leq \sum_{j\in J\setminus I} \Vert T e_j \Vert_{\mathcal H}.$$
\end{lemma}
\noindent {\it Proof}. By the min-max principle for bounded self-adjoint operators, we have
$$\mu_{n+1}(T)=\min_{\mathrm{dim}(F)=n} \max_{w \in F^{\perp},\Vert w \Vert=1} \langle \sqrt{T^*T}w,w\rangle.$$
Set $F=\mathrm{Span}\{ e_j,\ j\in I\}$. Given $w=\sum_{j\not \in I} c_j e_j$ with $\sum_j \vert c_j \vert^2=1$, we have by Schwarz inequality
$$\vert \langle \sqrt{T^*T}w,w\rangle\vert \leq \Vert \sqrt{T^*T}(w)\Vert=\Vert T(w)\Vert\leq \sum_{j\not \in I} \Vert T(e_j)\Vert,$$
and the proof is done. $\square$ 

\bigskip
Clearly since we have for all $n$ (we use Proposition \ref{conv}),
$$\mu_n(\lt_q)\leq \Vert \lt_q \Vert_{r,M} \leq L \Vert e^{2i\pi q \tau} \Vert_{H^2_{\widetilde{r}}} \Vert \lt_0 \Vert_{r,M}\leq 
L e^{2\pi \vert q\vert \Vert \tau \Vert_{r',\infty}} \Vert \lt_0 \Vert_{r,M},$$
for $r'>\widetilde{r}>K(M)r$, it is enough to estimate the singular values for all $n$ large enough. 
Let us set for all $R\geq 1$,
$$N_\infty(R):=\#\{ \alpha \in \Z^2\ :\ \Vert \alpha \Vert_{\infty} \leq R \},$$
where $\Vert \alpha\Vert_\infty=\max\{\vert \alpha_1\vert,\vert \alpha_2 \vert \}$.
We do have for all $R\geq 1$,
$$N_\infty(R)\leq N_\infty( [R]+1)=(2[R]+3)^2\leq 19R^2,$$
so that for all $n\geq19$, 
$$n\geq N_\infty \left(\sqrt{\frac{n}{19}} \right).$$
Therefore, by Lemma \ref{singular}, we have for all $n\geq 19$,
$$\mu_{n+1}(\lt_q)\leq \mu_{N_\infty \left(\sqrt{\frac{n}{19}} \right)+1} \leq \sum_{\Vert \alpha \Vert_{\infty}>\sqrt{\frac{n}{19}}}
\Vert \lt_q(\rho_\alpha)\Vert_{r,M}.$$ 
By Proposition \ref{conv}, we can write 
$$\Vert \lt_q(\rho_\alpha)\Vert_{r,M}\leq L e^{2\pi \vert q\vert \Vert \tau \Vert_{r',\infty}} \Vert \lt_0(\rho_\alpha) \Vert_{r,M}.$$
Using Theorem \ref{Adam}, we have 
$$\Vert \lt_0(\rho_\alpha) \Vert_{r,M}^2\leq \sum_{\beta \in \Z^2} \vert \langle \lt_0(\rho_\alpha),\rho_\beta \rangle_{r,M}\vert^2$$
$$\leq \sum_{\beta} e^{-2r_0(\vert \alpha\vert+\vert \beta \vert)}\leq C(r_0) e^{-2r_0 \vert \alpha \vert}.$$
As a consequence, we get
$$\mu_{n+1}(\lt_q)\leq C  e^{2\pi \vert q\vert \Vert \tau \Vert_{r',\infty}}\sum_{\Vert \alpha \Vert_{\infty}>\sqrt{\frac{n}{19}}}
e^{-r_0\vert \alpha \vert}.$$
On the other hand, we have by Stieltjes integration by parts (we also use the fact that $\vert \alpha\vert \geq \Vert \alpha \Vert_{\infty}$)
$$\sum_{\Vert \alpha \Vert_{\infty}>\sqrt{\frac{n}{19}}}
e^{-r_0\vert \alpha \vert}=\int_{\sqrt{\frac{n}{19}}}^\infty e^{-r_0 u}dN_\infty(u)=O\left( ne^{-r_0\sqrt{\frac{n}{19}} } \right)
=O\left (e^{-\kappa \sqrt{n}} \right),$$
for some $\kappa>0$, and all $n$ large. We have now obtained, at the cost of a large constant $C>0$, that for all $n,q$,
$$\mu_{n}(\lt_q)\leq C  e^{2\pi \vert q\vert \Vert \tau \Vert_{r',\infty}} e^{-\kappa \sqrt{n}}.$$
We now use a Weyl inequality (see \cite{Bsimon}, Thm 1.14) to transfer this singular value estimate into an actual eigenvalue
estimate. Indeed for all $N\geq 1$, we have
$$\vert \lambda_N(\lt_q)\vert^N\leq \prod_{n=1}^N \vert \lambda_n(\lt_q)\vert\leq \prod_{n=1}^N\mu_n(\lt_q).$$
This yields immediately
$$\vert \lambda_N(\lt_q)\vert\leq C  e^{2\pi \vert q\vert \Vert \tau \Vert_{r',\infty}} e^{-\frac{\kappa}{N}\sum_{n=1}^N\sqrt{n}}.$$
Since we obviously have
$$\sum_{n=1}^N \sqrt{n}\geq \int_0^N \sqrt{u}du=\frac{2}{3}N^{3/2},$$
the claim is proved.

\subsection{Trace formulas and spectral radius}
In this $\S$, we prove claim $(5)$ and then claim $(2)$ of Proposition \ref{blackbox}. First we start with the trace formula whose proof
runs exactly as in \cite{AlexAdam, FaureRoy},  and therefore we give only the outline.
By the general theory of trace class operators, see \cite{Bsimon}, chapter 3, we have
$$\mathrm{Tr}(\lt^n_q)=\lim_{N\rightarrow +\infty} \sum_{\Vert \alpha \Vert_\infty\leq N} \langle\lt^n(\rho_\alpha),\rho_\alpha \rangle_{r,M}.$$
On the other hand, 
$$\langle\lt^n(\rho_\alpha),\rho_\alpha \rangle_{r,M}=\langle A_{r,M} \lt^n(\rho_\alpha), A_{r,M} \rho_\alpha \rangle_{L^2(\T)}$$
$$=\int_{\T} e^{2i\pi q \tau^{(n)}(x)} \left(e_\alpha\circ A^n\right)(x) e_{-\alpha}(x)dm(x).$$
Therefore,
$$ \sum_{\Vert \alpha \Vert_\infty\leq N} \langle\lt^n(\rho_\alpha),\rho_\alpha \rangle_{r,M}=\int_{\T} e^{2i\pi q \tau^{(n)}(x)}D_N(A^nx-x)dm(x),$$
where 
$$D_N(x)=\left( \frac{\sin((2N+1)\pi x_1)}{\sin(\pi x_1)}\right)\left( \frac{\sin((2N+1)\pi x_2)}{\sin(\pi x_2)}\right)$$
is the $2$-dimensional Dirichlet kernel. This kernel tends (in distributional sense) to the Dirac measure at $0$ as $N\rightarrow +\infty$. Combining it with the fact that $A^n-Id$ is a local diffeomorphism, one can use a smooth partition of unity and a local change of coordinates to obtain as $N\rightarrow+\infty$
$$\mathrm{Tr}( \lt_q^n)=\sum_{A^nx=x} \frac{e^{2i\pi q \tau^{(n)}(x)}}{\vert \det(I-D_xA^n)\vert}.$$
Notice that these spectral identities involve only sums over periodic points and do not depend on the choice of the space $\Aniso$, and hence the spectrum. We can now prove claim $(2)$. While some tedious calculation involving the norms $\Vert \lt_q^n\vert_{r,m}$ can definitely lead to the fact that for all $q$, $\rho(\lt_q)\leq 1$, we choose an easier route here. Because $\lt_q$ are trace class operators, we can consider the
Fredholm determinants
$$\mathcal{Z}_q(\zeta):=\det(I-\zeta \lt_q).$$
Each determinant $\mathcal{Z}_q(\zeta)$ is an entire function of $\zeta \in \C$, whose zeros are given by the inverses of the 
non-zero {\it eigenvalues} of $\lt_q$. For all $\vert \zeta \vert$ small enough, we have the identity (by the trace formula and Lidskii theorem)
$$\mathcal{Z}_q(\zeta)=\exp \left ( \sum_{n=1}^\infty \frac{\zeta^n}{n} \sum_{A^nx=x} \frac{e^{2i\pi q \tau^{(n)}(x)}}{\vert \det(I-D_xA^n)\vert}  \right).$$
However by Lemma \ref{pressure1}, we know that for all $\epsilon>0$ and $n$ large,
$$\left \vert \sum_{A^nx=x} \frac{e^{2i\pi q \tau^{(n)}(x)}}{\vert \det(I-D_xA^n)\vert} \right \vert \leq C {e^{n(P(-\log J^u)+\epsilon)}}.$$
But by Bowen \cite{BowenBook}, chapter 4, we know that $P(-\log J^u)=0$, therefore the series
$$\sum_{n=1}^\infty \frac{\zeta^n}{n} \sum_{A^nx=x} \frac{e^{2i\pi q \tau^{(n)}(x)}}{\vert \det(I-D_xA^n)\vert}$$
are absolutely convergent for all $\vert \zeta \vert <1$, which proves by analytic continuation that $\mathcal{Z}_q(\zeta)$ cannot vanish
if $\vert \zeta\vert <1$. This is enough to conclude that for all $q$, all the eigenvalues of $\lt_q$ have modulus smaller or equal to $1$, hence proving claim $(2)$.
\subsection{Lebesgue measure on ideals of $\Aniso$}
Below we prove claim $(4)$ of Proposition \ref{blackbox}.

\bigskip \noindent
For all trigonometric polynomial $\psi$, we have 
$$L_m(\psi):=\int_{\T}\psi dm=\langle A_{r,M}\psi,A_{r,M} e_0 \rangle_{L^2(\T)}=\langle \psi, \rho_0 \rangle_{r,M}.$$
This shows that one can clearly extend the functional $L_m$ by density to $\Aniso$ by setting for all $G\in \Aniso$,
$$L_m(G):=\langle G, \rho_0 \rangle_{r,M}.$$
Let $g=\sum_\alpha C_\alpha \rho_\alpha$ be a trigonometric polynomial and 
$\varphi=\sum_{\alpha}B_\alpha \rho_\alpha \in \Aniso\setminus\{0\}$. We can check that we have the identity
$$L_m(\varphi g)=\langle \varphi g, \rho_0 \rangle_{r,M}=\sum_{\alpha \in \Z^2} C_\alpha B_{-\alpha} e^{4\pi r (\vert \alpha^+\vert-\vert \alpha^-\vert)}.$$
Because $\varphi\neq 0$, 
at least one of the coefficients $B_\beta$ is non vanishing, say $B_{\beta_0}$. Now simply pick $g=\rho_{-\beta_0}$
so that all $C_\alpha$ are vanishing except $C_{-\beta_0}=1$. We have 
$$L_m(\varphi g)=B_{\beta_0}e^{4\pi(\vert \beta_0^+\vert -\vert \beta_0^-\vert)}\neq 0,$$
and the proof is done. $\square$

We would like to conclude by an open question related to the spaces $\Aniso$. If $A$ is not volume preserving, then the SRB-measure still makes sense as a continuous functional
$L_{SRB}$ acting on $\Aniso$. In that case it is no longer absolutely continuous but is a radon measure with full support (consequence of topological mixing). 
Can one show that given $\varphi \neq 0$ in $\Aniso$, there exists a trigonometric
polynomial $g$ such that $L_{SRB}(\varphi g)\neq 0$ ? The problem is likely to boil down to an analytic wave front set argument and checking that the (distributional) pairing $L_{SRB}\otimes \varphi$
is not zero.


\begin{thebibliography}{10}

\bibitem{AlexAdam}
Alexander Adam.
\newblock Generic non-trivial resonances for anosov diffeomorphism.
\newblock {\em Preprint 2016}.

\bibitem{BaladiBook1}
Viviane Baladi.
\newblock {\em Positive transfer operators and decay of correlations},
  volume~16 of {\em Advanced Series in Nonlinear Dynamics}.
\newblock World Scientific Publishing Co., Inc., River Edge, NJ, 2000.

\bibitem{BaladiBook2}
Viviane Baladi.
\newblock {\em Dynamical zeta functions and dynamical determinants for
  hyperbolic maps, a functional approach}.
\newblock Springer, 2016.

\bibitem{BowenBook}
Rufus Bowen.
\newblock {\em Equilibrium states and the ergodic theory of {A}nosov
  diffeomorphisms}, volume 470 of {\em Lecture Notes in Mathematics}.
\newblock Springer-Verlag, Berlin, revised edition, 2008.
\newblock With a preface by David Ruelle, Edited by Jean-Ren{\'e} Chazottes.

\bibitem{Brin1}
M.~I. Brin.
\newblock The topology of group extensions of {$C$}-systems.
\newblock {\em Mat. Zametki}, 18(3):453--465, 1975.

\bibitem{BurnsWilkinson}
Keith Burns and Amie Wilkinson.
\newblock Stable ergodicity of skew products.
\newblock {\em Ann. Sci. \'Ecole Norm. Sup. (4)}, 32(6):859--889, 1999.

\bibitem{BW1}
Oliver Butterley and Khadim War.
\newblock Open sets of exponentially mixing anosov flows.
\newblock {\em Preprint 2016}.

\bibitem{Dolgopyat}
Dmitry Dolgopyat.
\newblock On mixing properties of compact group extensions of hyperbolic
  systems.
\newblock {\em Israel J. Math.}, 130:157--205, 2002.

\bibitem{FaureRoy}
Fr{\'e}d{\'e}ric Faure and Nicolas Roy.
\newblock Ruelle-{P}ollicott resonances for real analytic hyperbolic maps.
\newblock {\em Nonlinearity}, 19(6):1233--1252, 2006.

\bibitem{Liverani}
Carlangelo Liverani.
\newblock On contact {A}nosov flows.
\newblock {\em Ann. of Math. (2)}, 159(3):1275--1312, 2004.

\bibitem{Naud2}
Fr{\'e}d{\'e}ric Naud.
\newblock Analytic continuation of a dynamical zeta function under a
  {D}iophantine condition.
\newblock {\em Nonlinearity}, 14(5):995--1009, 2001.

\bibitem{NaudHP}
Fr{\'e}d{\'e}ric Naud.
\newblock Entropy and decay of correlations for real analytic semi-flows.
\newblock {\em Ann. Henri Poincar\'e}, 10(3):429--451, 2009.

\bibitem{PalisTakens}
Jacob Palis and Floris Takens.
\newblock {\em Hyperbolicity and sensitive chaotic dynamics at homoclinic
  bifurcations}, volume~35 of {\em Cambridge Studies in Advanced Mathematics}.
\newblock Cambridge University Press, Cambridge, 1993.
\newblock Fractal dimensions and infinitely many attractors.

\bibitem{ParryPollicott}
William Parry and Mark Pollicott.
\newblock Zeta functions and the periodic orbit structure of hyperbolic
  dynamics.
\newblock {\em Ast\'erisque}, (187-188):268, 1990.

\bibitem{Bsimon}
Barry Simon.
\newblock {\em Trace ideals and their applications}, volume 120 of {\em
  Mathematical Surveys and Monographs}.
\newblock American Mathematical Society, Providence, RI, second edition, 2005.

\bibitem{BJS2}
Julia Slipantschuk, Oscar~F. Bandtlow, and Wolfram Just.
\newblock Complete spectral data for analytic anosov maps of the torus.
\newblock {\em Preprint 2016}.

\bibitem{Tsujii}
Masato Tsujii.
\newblock Exponential mixing for generic volume-preserving anosov flows in
  dimension three.
\newblock {\em Preprint}, 2016.

\end{thebibliography}
\end{document}